\documentclass[11pt,a4paper]{article}
\usepackage{a4wide}
\usepackage{amsfonts,amsmath,amsthm,amssymb}
\usepackage{amstext}
\usepackage{graphicx,geompsfi-pdf,subfigure}
\usepackage[small,nohug,heads=littlevee]{diagrams} 
\diagramstyle[labelstyle=\scriptstyle] 
\newcommand{\e}{\varepsilon}
\let\epsilon\varepsilon

\makeatletter
\@addtoreset{equation}{section}
\makeatother
\newtheorem{lemma}{Lemma}[section]
\newtheorem{theorem}[lemma]{Theorem}

\newtheorem{conjecture}[lemma]{Conjecture}

\newcommand{\weak}{\rightharpoonup}

\newcommand{\beq}{\begin{equation}}
\newcommand{\eeq}{\end{equation}}

\def\k{\mathbf k}
\def\0{\mathbf 0}

\long\def\drop#1{}
\def\pref#1{(\ref{#1})}
\DeclareMathOperator{\supp}{supp}

\def\R{\mathbf R}
\def\N{\mathbf N}
\def\Z{\mathbf Z}
\def\H{\mathcal H}
\def\Hmo#1{\|#1-{\textstyle\dashint #1 }\|_{H^{-1}}}
\def\Hmospace#1#2{\|#1-{\textstyle\dashint #1 }\|_{H^{-1}(#2)}}

\def\Ees{E_{\e,\sigma}}

\def\Eee{E_{\e,\eta}}
\def\Eeen{E_{\e_n,\eta_n}}
\def\Eeet{E^{\rm 2D}_{\e,\eta}}
\def\Eeent{E^{\rm 2D}_{\e_n,\eta_n}}

\def\sfE{\mathsf{E}}
\def\Fea{{{\sf{E_\eta^{\mathrm{2D}}}}}}
\def\Feb{{\sf{E_\eta}}}
\def\Fen{{\sfE_{\eta_n}}}
\def\Fza{{\sf{E_0^{\mathrm{2D}}}}}
\def\Fzb{{\sf{E_0}}}

\def\fza{{e_0^{\rm 2D}}}
\def\fzb{e_0}
\def\lscfz{\overline e_0}
\def\lscfza{\overline{e_0^{\rm 2D}}}
\def\sfF{\mathsf{F}}
\def\Hea{{\sf{F_\eta^{\mathrm{2D}}}}}
\def\Heb{{\sf{F_\eta}}}
\def\Hen{{\sfF_{\eta_n}}}
\def\Hza{{\sf{F_0^{\mathrm{2D}}}}}
\def\Hzb{{\sf{F_0}}}
\def\Jee{F_{\e,\eta}}
\def\Jeen{F_{\e_n,\eta_n}}
\def\Jeet{F^{\rm 2D}_{\e,\eta}}
\def\Jeent{F^{\rm 2D}_{\e_n,\eta_n}}
\def\T{\mathbf{T}}
\def\Ta{{\mathbf{T}^2}}
\def\Tb{{\mathbf{T}^3}}
\def\logeta{\left|\log\eta\right|}
\def\invlogeta{\left|\log\eta\right|^{-1}}

\def\Xint#1{\mathchoice
   {\XXint\displaystyle\textstyle{#1}}%
   {\XXint\textstyle\scriptstyle{#1}}%
   {\XXint\scriptstyle\scriptscriptstyle{#1}}%
   {\XXint\scriptscriptstyle\scriptscriptstyle{#1}}%
   \!\int}
\def\XXint#1#2#3{{\setbox0=\hbox{$#1{#2#3}{\int}$}
     \vcenter{\hbox{$#2#3$}}\kern-.5\wd0}}

\def\dashint{\Xint-}

\newenvironment{remark}%
  {\par\medbreak\refstepcounter{lemma}%
    \noindent\textbf{Remark~\thetheorem.}}%
  {\par\medskip}

\begin{document}

\begin{title}{Small Volume Fraction Limit of the Diblock Copolymer Problem: II. Diffuse-Interface Functional}

\author{Rustum Choksi\footnote{Department of Mathematics, Simon Fraser University, 
Burnaby, Canada, choksi@math.sfu.ca}  \and Mark A. Peletier\footnote{Department of Mathematics and Institute for Complex Molecular Systems, Technische Universiteit Eindhoven, The Netherlands, m.a.peletier@tue.nl}}
\end{title}

\maketitle 

\begin{abstract} 
We present the second of two articles on the small volume fraction limit of a nonlocal Cahn-Hilliard functional introduced  to model microphase separation of diblock copolymers. 
After having established the results for the sharp-interface version of the functional (\cite{ChoksiPeletier09}), we consider here the full diffuse-interface functional   and address the 
 limit in which $\epsilon$ and  the volume fraction tend to zero but the number of minority phases (called {\it particles}) remains $O(1)$. 
Using the language of $\Gamma$-convergence, we focus on   two levels of this convergence, and derive first- and second-order {\it effective}  energies, whose energy landscapes are simpler and more transparent. These limiting energies are only finite on weighted sums of delta functions, corresponding to the concentration of mass into `point particles'.
At the highest 
level, the effective energy is entirely local and contains information about the size of each particle but no information about their spatial distribution. 
 At the next level we encounter a Coulomb-like 
interaction between the particles, which is responsible for the pattern formation. 
We present the results in  three  dimensions and comment on their two-dimensional analogues.


\medskip


\end{abstract}

\tableofcontents
\section{Introduction}
\subsection{The Functional}
This paper is concerned with asymptotic properties of  
the following  nonlocal Cahn-Hilliard  energy functional defined on $H^1(\R^d)$: 
\begin{equation}
\label{def:Ees}
{\cal E}(u) \;:=\; \e \,  \int_{\T^d} \,  |\nabla u|^2  \, d { x} \, \, + \,\, \frac1\e\int_{\T^d} W(u) \, d {x}
 \, \, + \,\, \gamma 
 \, \Hmospace{u}{\T^d}^2, 
\end{equation} 
where we take the  double-well potential $W (u) := u^2(1-u^2)$. 
Here the order parameter $u$ is defined on the flat torus $\T^d=\R^d/\Z^d$, i.e.\ the square $[-\frac{1}{2},\frac{1}{2}]^d$ with periodic boundary conditions, and has two preferred states $u = 0$ and $u = 1$. We are interested in the structure of minimizers of  over $u$ with fixed mass $\dashint_{\T^d} u = {  f}$ where $f \in (0,1)$. 
The first term $\e\int |\nabla u|^2$ penalizes large gradients, and acts as a counterbalance to the second term, smoothing the `interface' that separates the two phases. The third (nonlocal) term  is defined as
\[
\Hmospace{u}{\T^d}^2 = \int_{\T^d} |\nabla w|^2 \, dx,
\qquad\text{where}\qquad
-\Delta w = u - \dashint_{\T^d} u  .
\]
This term favors high-frequency oscillation, as can be recognized in the $1/|\k|^2$-penalization in a Fourier representation:
\[
\Hmospace{u}{\T^d}^2 = \sum_{\k\in \Z^d\setminus\{0\}} 
  \frac {|\hat u(\k)|^2}{4\pi^2 |\k|^2}.
\]
If the parameter $\gamma$ is large enough, this term may push the system away from large, bulky structures, and favor variation and oscillation at intermediate scales, i.e. give rise to patterns with an intrinsic length scale. As we explain in the sequel, we refer to this mass-constrained variational problem as {\it the diblock copolymer problem}. 
When the mass constraint 
$f$ is close to $0$ or $1$, minimizing patterns will consist of small inclusions of one phase in a large `sea' of the other. 
We wish to explore this regime  via 
the asymptotic behavior of the functional in a limit wherein 
\begin{itemize}
\item both $\epsilon$ and the volume/mass fraction $f$ of the minority phase tend to zero (appropriately slaved together)
\item $\gamma$ is chosen in order to keep the number of minority phase {\it particles}  $O(1)$. 
\end{itemize}
We will primarily concern ourselves with the case $d=3$ but remark on the analogous results for  $d = 2$.

\subsection{The Spherical Phase in Diblock Copolymers}
The functional ${\cal E}$ was introduced by Ohta and Kawasaki to model self-assembly of diblock copolymers 
\cite{OhtaKawasaki86, NishiuraOhnishi95}. The nonlocal term  is associated with long-range interactions and  
connectivity of the sub-chains in the diblock copolymer  
macromolecule\footnote{See \cite{ChoksiRen03} for a derivation and  the relationship to the physical material parameters and basic models for inhomogeneous polymers. Usually the wells are taken to be $\pm 1$ representing pure phases of $A$ and $B$-rich regions. For convenience, we have rescaled to wells at $0$ and $1$.}. 
The order parameter $u$ represents the relative monomer density, with $u = 0$ corresponding to a pure-$A$ region and $u = 1$ to a pure-$B$ region. The interpretation of $f$ is therefore the relative abundance of the $A$-parts of the molecules, or equivalently the volume fraction of the $A$-region. The constraint of fixed average $f$ reflects that in an experiment the composition of the molecules is part of the preparation and does not change during the course of the experiment. 
From \pref{def:Ees} the incentive for pattern formation is clear: the first term penalizes oscillation, the second term favors separation into regions of $u=0$ and $u=1$, and the third favors rapid oscillation. Under the mass constraint the three can not vanish simultaneously, and the net effect is to set a fine scale structure depending on $\e, \gamma$ and $f$.
The precise geometry of the phase separation (i.e.\ the information contained in a minimizer of~\pref{def:Ees}) depends largely on the volume fraction~$f$. In fact, as explained in \cite{ChoksiPeletierWilliams09}, the two natural parameters controlling the phase diagram are $\Gamma = (\e^{3/2} \sqrt{\gamma})^{-1}$ and $f$. When $\Gamma$ is large and $f$ is close to $0$ or $1$, numerical experiments~\cite{ChoksiPeletierWilliams09} and experimental observations~\cite{BatesFredrickson99} 
reveal structures resembling {\it small well-separated spherical regions of the minority phase}. We often refer to  such small regions as {\it particles}, and they are the central objects of study of this paper. 
Since we are interested in a regime of small volume fraction, it seems natural to seek asymptotic results.  Building on our previous work in \cite{ChoksiPeletier09}, it is the purpose of this article to give a rigorous 
asymptotic description of the energy in a limit wherein the volume fraction tends to zero but where the number of particles in a minimizer remains $O(1)$. That is, we examine the limit where minimizers converge to weighted Dirac delta point measures and seek effective energetic descriptions for their positioning and local structure.

\begin{figure} 
\centerline{{\psfig{height=2in,figure=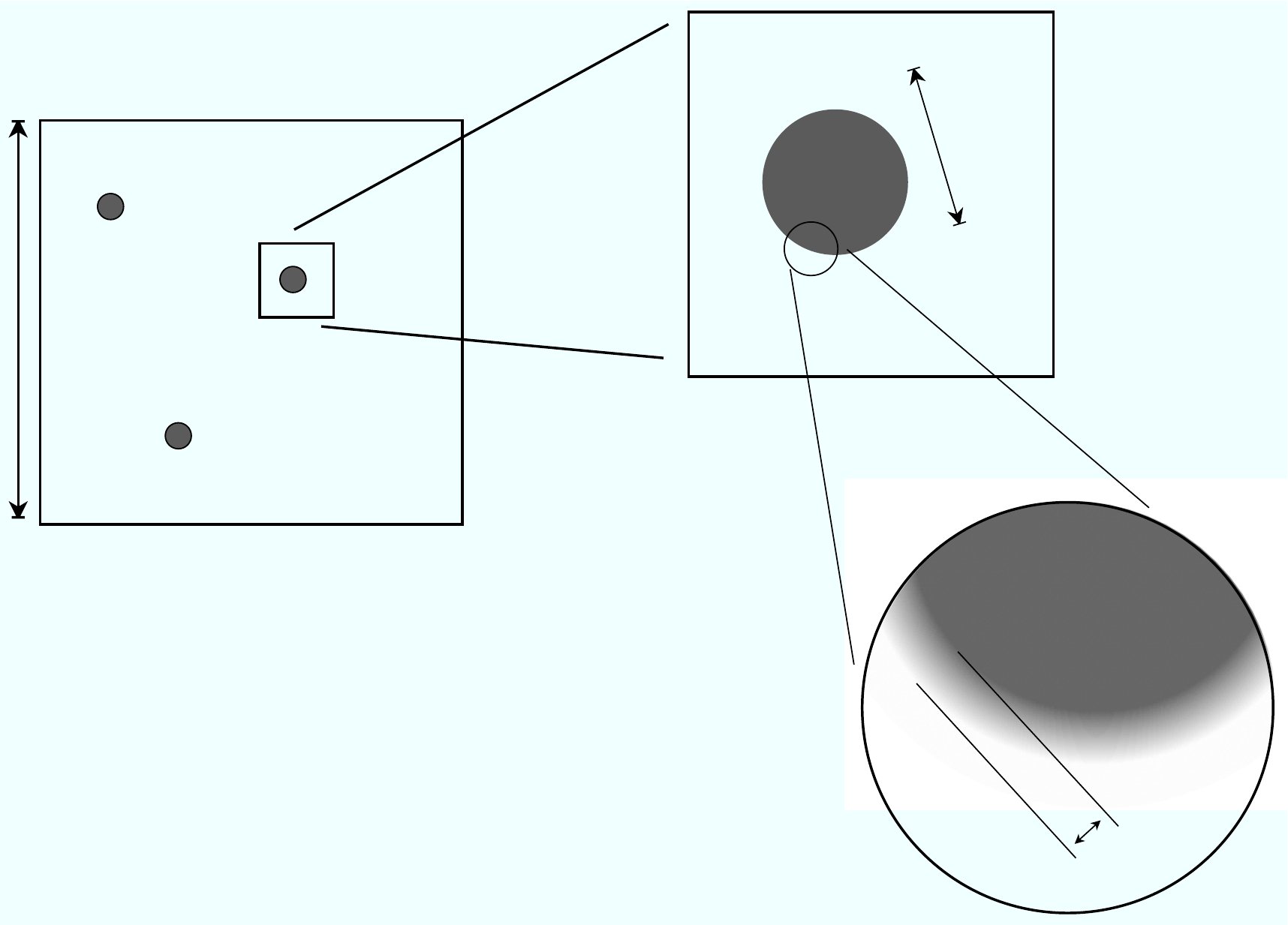}}}
\caption{A two-dimensional cartoon of small-particle structures}
\label{fig:scales}
\end{figure}

The small particle structures of this paper are illustrated (for two space dimensions) in  Figure~\ref{fig:scales}. 
There are three length scales involved: the large scale of the periodic box $\T^d$, the intermediate scale of the droplets, and the smallest scale of the thickness of the interface. Two of these scales are known beforehand: the the size of the box we have chosen to be $1$, and the interfacial thickness should be $O(\e)$ by the discussion above. The intermediate scale $\ell$, the size of the droplets, is not yet fixed, and will depend on the two remaining parameters: the parameter~$\gamma$ in ${\cal E}$ and the volume fraction $f$. 


For a function $u$, the \emph{mass} is defined as $f = \int_{\T^d} u$. In Fig.~\ref{fig:scales} the region where $u\approx 1$ is small, suggesting that $\int_{\T^d} u$ is small. We characterize this by introducing a parameter $\eta$, which will tend to zero, and by assuming that the mass $\int_{\T^d} u$ tends to zero at the rate of $\eta^d$: 
\begin{equation}
\label{def:eta}
f \, = \, \int_{\T^d} u =  M\eta^d, \qquad \text{for some fixed }M>0.
\end{equation}
After rescaling with respect to $\eta$, $M$ will be the mass of the rescaled functions. We now have three parameters $\e$, $\gamma$, and $\eta$, which together determine the behavior of structures under the energy $\Ees$; in this paper we keep the parameter $M$ fixed. Let us fix $d = 3$. 
In Section \ref{sect-rescale} we see that in terms of $v: = u / \eta^3$, the relevant functional is 
\[ \Eee (v) \, : = \,  \eta \left[ \e \, \eta^3 \int_{\T^3} \,  |\nabla v|^2  \, d { x} \, \, + \,\, \frac{\eta^3}{\e}\int_{\T^3} \widetilde W (v) \, d {x} \right]  \,\, + \,\, \eta  \Hmospace v\Tb ^2,\]
where $\widetilde W (v) := v^2 ( 1 - \eta^3 v)^2$.
Via a suitable slaving of $\epsilon$ to $\eta$ (see Theorem \ref{T1}), we prove, via $\Gamma$-convergence,  a rigorous asymptotic expansion for $E_{\e(\eta), \eta} $ of the form: 
\[ E_{\e(\eta), \eta} \, = \, \Fzb  \,\, + \,\, \eta \Hzb \,\, + \,\, \text{higher order terms}, \]
where both $\Fzb$ and $\Hzb$ are defined over weighted Dirac point masses and may be viewed as {\it effective} energies at the first and second order. Their essential properties can be summarized as follows:   
\begin{itemize}
\item
$\Fzb$, the effective energy at the highest level,  is entirely local: it is the sum of local energies of each particle, and is blind to the spatial distribution of the particles. 
The particle effective energy only depends on the mass of that particle. 

\item $\Hzb$, the effective energy at the next level, contains a Coulomb-like 
interaction between the particles. It is this latter part of the energy which we expect 
enforces a periodic array of particles. 
\end{itemize} 

The proof a Theorem \ref{T1} relies heavily on our previous work for the sharp interface limiting functional  $\Feb$ (see Section 
\ref{sharp} for its precise definition) obtained by fixing $\eta$  in $\Eee$ and letting $\e $ tend to zero. 
The well-known Modica-Mortola Theorem~\cite{ModicaMortola77}  makes this limit $\Feb$ 
precise in the sense of $\Gamma$-{convergence}. 
The small-$\eta$ asymptotics of $\Feb$  were proved in \cite{ChoksiPeletier09}, and 
the main result of this article (Theorem \ref{T1}) is to establish  the same limiting behavior but in the diagonal limit of both $\epsilon$ and $\eta$ tending to zero.  We summarize these limits (for the leading order)  in the diagram below. 
 \begin{diagram} 
\Eee&\\ 
\dTo^{ \substack{\eta \text{ fixed},\,\,  \epsilon \rightarrow 0:\\\text { Modica-Mortola Theorem}} }
&\rdTo^{\epsilon (\eta) \, \& \,  \eta \,  \longrightarrow \,  0: \,\, \text{Theorem } \ref{T1}} & \\ 
\Feb&\rTo_{\substack{\eta \, \longrightarrow \, 0:\\  \,\, \text{ Theorem  } \ref{T2}   \text {, proved in  \cite{ChoksiPeletier09}} } \qquad  \qquad } &\Fzb 
\end{diagram}  
The article is organized as follows. In Section~\ref{sect-rescale}, we discuss the rescalings and state the main result Theorem \ref{T1}. 
Section~\ref{sharp} explicitly states the main results of our previous paper  \cite{ChoksiPeletier09} which form the basis for the proof of Theorem \ref{T1} presented in Section~\ref{prooftheorem}.  In Section~\ref{localstructure}, we discuss the variational problem associated with the first order $\Gamma-$limit $\Fzb$, connecting it with an old problem of Poincar\'{e} and presenting some conjectures.  In Section~\ref {2D}, we discuss the necessary modifications in two dimensions.


\bigskip

\section{Some definitions and notation}
\label{sec:notation}

We recall the definitions and notation of \cite{ChoksiPeletier09}. 
We use $\T^{d}=\R^d/\Z^d$ to denote the $d$-dimensional flat torus of unit volume. We will primarily be concerned with the case $d = 3$. 
For the use of convolution we note that $\T^d$ is an additive group, with neutral element $0\in\T^d$ (the `origin' of $\T^d$).
  For $v \in BV (\T^d; \{0,1\})$ we denote by 
\[ \int_{\T^d} |\nabla v|\]
the total variation measure evaluated on $\T^d$, i.e. $\|\nabla u\| (\T^d)$ (see e.g.~\cite{Ambrosio89} or~\cite[Ch.~3]{AmbrosioFuscoPallara00}). 
Since $v$ is the characteristic function of some set $A$, it is simply 
the notion of its perimeter.  
Let $X$ denote the space of Radon measures on $\T^d$. For $ \mu_\eta, \mu \in X$, 
$ \mu_\eta \weak \mu$ denotes weak-$\ast$ measure convergence, i.e.
\[ 
\int_{\T^d} f \, d \mu_\eta \, \, \rightarrow \, \,   \int_{\T^d} f \, d \mu
\]
for all  $f \in C(\T^n)$.  We use the same notation for functions, i.e. 
when writing  $v_\eta\weak v_0$, we interpret $v_\eta$ and $v_0$ as measures whenever necessary.

We  introduce the Green's function  $G_{\T^d}$ for $-\Delta$ in dimension $d$ on $\T^d$. It is the solution of 
\[
- \Delta G_{\T^d} = \delta \, - \, 1,
\qquad\text{with}\qquad
\int_{\T^d} G_{\T^d} = 0,
\]
where $\delta$ is the Dirac delta function at the origin. 
In three dimensions\footnote{In two dimensions, the Green's function $G_{\T^2}$ satisfies
\begin{equation}
\label{eq:G_T-g2}
G_{\T^2} (x) \, = \, -\frac1{2\pi} \log |x| \, + \, { g}^{(2)}(x)
\end{equation}
for all $x=(x_1,x_2)\in \R^2 $ with $\max\{|x_1|,|x_2|\}\leq 1/2$, where the function $g^{(2)}$ is continuous on $[-1/2,1/2]^2$ and  $C^\infty$ in a neighborhood of 
the origin.}, we have 
 \begin{equation}
\label{eq:G_T-g3}
G_{\T^3} (x) \, = \, \frac1{4 \pi |x|}\,  +\,  g^{(3)}(x)
\end{equation}
for all $x=(x_1,x_2, x_3)\in \R^3 $ with $\max\{|x_1|,|x_2|, |x_3|\}\leq 1/2$, where the function $g^{(3)}$  is continuous on $[-1/2,1/2]^3$ and smooth in a neighborhood of the origin.

For $\mu \in X$ such that $\mu (\T^d) = 0$, we may solve 
\[ 
-\Delta w \, = \, \mu,  
\]
in the sense of distributions on $\T^d$.
If $w \in H^1 (\T^d)$, then  $\mu \in H^{-1} (\T^d)$, and 
 \[ \|\mu\|_{H^{-1}(\T^d)}^2 \, := \, \int_{\T^d} |\nabla w|^2 \, dx.  \]
In particular, if $u \in L^2({\T^d})$ then $\left(u - \dashint u\right) \in H^{-1} (\T^d)$ and   
\[\Hmospace{u}{\T^d}^2 
 \, = \, \int_{\T^d} \int_{\T^d} u(x)u(y) \, G_{\T^d} (x-y)\, dx\, dy.
\]
Note that on the right-hand side we may write the function $u$ rather than its zero-average version $u-\dashint u$, since the function $G_{\T^d}$ itself is chosen to have zero average. 

If $f$ is the characteristic function of a set of finite perimeter on all of $\R^3$, we define 
\[ \|f\|_{H^{-1}({\R^3})}^2 \, = \, \int_{\R^3}\int_{\R^3} \frac{f(x) \, f(y)}{4 \pi |x-y|}\,   dx\, dy.  \] 

\section{Rescalings and Statements of the Results}
\label{sect-rescale}
We now rescale the energy $\mathcal E$ in \eqref{def:Ees}. 
Starting in three dimensions, for $\eta >0$ we define 
\[ v\, :=\, \frac{u}{\eta^3}, \]
so that $\mathcal E$ becomes in terms of $v$
\beq\label{rescaledenergy}
\; \e \, \eta^6 \int_{\T^3} \,  |\nabla v|^2  \, d { x} \, \, + \,\, \frac{\eta^6}{\e}\int_{\T^3} \widetilde W (v) \, d {x}
 \, \, + \,\, \gamma \, \eta^6 \, 
 \, \Hmospace{v}{\T^3}^2, 
 \eeq
 where 
 \[ \widetilde W (v) := v^2 ( 1 - \eta^3 v)^2. \]

In order to find the correct scaling of $\gamma$ in terms of  $\eta$,  we consider a collection $v_\eta: {\T^3} \rightarrow \{0, 1/\eta^3\}$ of  components of the form 
\begin{equation}\label{formv}
v_\eta \, = \, \sum_i v_\eta^i, \qquad v_\eta^i \, = \, 
\frac{1}{\eta^n} \chi^{}_{A_i},
\end{equation}
where the $A_i$ are disjoint, connected subsets of $\T^3$. 
Then under the assumption that the number of $A_i$ remains $O(1)$, we find 
\[ \eta \left[ \e \, \eta^3 \int_{\T^3} \,  |\nabla v|^2  \, d { x} \, \, + \,\, \frac{\eta^3}{\e}\int_{\T^3} \widetilde W (v) \, d {x} \right] 
\,\,  \stackrel{{\e \ll \eta}}{\sim} \,\, \eta  \int_{\T^3} \,  |\nabla v| \, = \, O(1). \]
Here we are using the well-known Modica-Mortola convergence theorem \cite{ModicaMortola77, Braides02} linking the perimeter to the scaled Cahn-Hilliard terms.
A simple calculation (done in \cite{ChoksiPeletier09}) shows that the leading order of the $\| v_\eta - \dashint v_\eta \|^2_{H^{-1}(\T^3)}$ is $1 / \eta$, and that this leading contribution is from the {\it self-interactions}, i.e. 
$\| v^i_\eta - \dashint v^i_\eta \|^2_{H^{-1}(\T^3)}$ is $1 / \eta$. 
Thus balancing the third term in \pref{rescaledenergy} implies choosing $\gamma \sim 1 / \eta^3$. Hence we set 
\[ \gamma \, = \, \frac{1}{\eta^3}.\]
Choosing the proportionality constant equal to $1$ entails no loss of generality, since in the limit $\e\to0$ this constant can be scaled into the mass $M$ defined in~\eqref{def:eta}.

With this choice, one finds 
 \[ {\cal E} (u) \,  = \, \eta^2 \left\{ \eta \left[ \e \, \eta^3 \int_{\T^3} \,  |\nabla v|^2  \, d { x} \, \, + \,\, \frac{\eta^3}{\e}\int_{\T^3} \widetilde W (v) \, d {x} \right]  \,\, + \,\, \eta  \Hmospace v\Tb ^2  \right\}, \]
noting that the contents of the outer parentheses is $O(1)$ as $\eta \rightarrow 0$ with $\e \ll \eta$. 
Thus let us define the re-normalized energy 
\begin{equation}
\label{def:rescaling_energy}
\Eee (v) \, : = \,  \eta \left[ \e \, \eta^3 \int_{\T^3} \,  |\nabla v|^2  \, d { x} \, \, + \,\, \frac{\eta^3}{\e}\int_{\T^3} \widetilde W (v) \, d {x} \right]  \,\, + \,\, \eta  \Hmospace v\Tb ^2 . 
\end{equation}
We are interested in the small-$\eta$ behavior of $\Eee$ and  describe this behavior via functionals defined over Dirac point masses.  
Let us first introduce the remaining relevant functionals in our analysis. First we define the surface tension 
\beq\label{st}
 \sigma \, : = \, 2 \int_0^1 \sqrt{W(t)} \, dt. \eeq
For the leading order, we define 
\begin{align}
\fzb(m) &:= \inf\left\{\sigma  \int_{\R^3} |\nabla z| \,\,\,\,  + \, \,\,\, \|z\|_{H^{-1}(\R^3)}^2: 
   z\in BV(\R^3;\{0,1\}), \int_{\R^3} z = m \right\}. 
\label{def:fe3}
\end{align}
and  
\begin{align*}
\Fzb(v) := \begin{cases}
\sum_{i=1}^\infty \fzb(m^i) & \text{if } v = \sum_{i=1}^\infty m^i \delta_{x^i}, \,\,\,  \{x^i\} \text{ distinct, and }m^i\geq0\\
\infty & \text{otherwise}.
\end{cases}
\end{align*}
For the next order we  note that  
among all measures of mass $M$, the global minimum of $\Fzb$ is given by
\[
\min \left\{\Fzb(v): \int_{\Tb} v = M\right\} = \fzb(M).
\]
We will recover the next term in the expansion as the limit of $\Eee - \fzb$, appropriately rescaled, that is of the functional
\[ 
\Jee(v_\eta) := \eta^{-1}\left[ \Eee(v_\eta) - \fzb\left(\int_\Tb v_\eta\right)\right] .
\]
Its limiting behavior will be characterized by the functional 
\[
\Hzb (v) := \begin{cases}
\displaystyle
\sum_{i=1}^\infty g^{(3)}(0) \, (m^i)^2 \, +\, & \\
 \qquad \sum_{{i\not=j}} m^im^j\, G_\Tb (x^i-x^j)
  & \displaystyle\text{if } v = \sum_{i=1}^n m^i \delta_{x^i} \text{ with } \{x^i\} \text{ distinct,  } \{m^i\}\in{\cal M}\\
\infty & \text{otherwise},
\end{cases}
\]
where $g^{(3)}$ is defined in~\eqref{eq:G_T-g3} and
\begin{eqnarray*}
{\cal M}
& := & \bigg\{ \{m^i\}_{i\in\N}: m^i\geq 0,\,    \fzb(m^i) \,\, \text{admits a minimizer for each $i$,}  \\
& & \qquad  \qquad \qquad \qquad \qquad \qquad \qquad \qquad 
\left.  \text{and} \, \,\,   \sum_{i=1}^\infty \fzb(m^i) = \fzb\Bigl(\sum_{i=1}^\infty m^i\Bigr)  \right\}.
\end{eqnarray*}

\bigskip

In Theorem \ref{T1} we prove that 
\[
\Eee\, \stackrel{\Gamma}\longrightarrow \, \Fzb  \qquad 
{\rm and} \qquad 
\Jee \, \stackrel{\Gamma}\longrightarrow \, \Hzb. 
\]
Precisely, 
\begin{theorem}
\label{T1}
\begin{itemize} 
\item (Condition 1 -- the lower bound and compactness) 
Let $\e_n$ and $\eta_n$ be sequences tending to zero such that, for some $\zeta >0$,  $\e_n = o(\eta_n^{4 + \zeta})$. Let $v_n$ be a sequence such that the energy 
$\Eeen(v_n)$ is bounded. Then
(up to a subsequence) $v_n\weak v_0$, $\supp v_0$ is countable, and 
\begin{equation}
\label{lb:diffuse}
\liminf_{n\to\infty} \Eeen(v_n) \geq \Fzb(v_0).
\end{equation}
If in addition $\Jeen(v_n)$ is bounded and $\zeta\geq1$, then the limit $v_0$ is a global minimizer of $\Fzb$ under constrained mass, and 
\begin{equation}
\label{lb:diffuse-nextlevel}
\liminf_{n\to\infty} \Jeen(v_n) \geq \Hzb(v_0).
\end{equation}
\item (Condition 2 -- the upper bound)
There exist two continuous functions $C_1,C_2:[0,\infty)\to[0,\infty)$ with $C_1(0)=C_2(0)=0$ with the following property.
Let $\e_n$ and $\eta_n$ be sequences tending to zero and let  $\e_n \leq C_1(\eta_n)$.  Let $v_0$ be such that $\Fzb(v_0)<\infty$.
Then there exists a sequence $v_n \weak v_0$ such that 
\begin{equation}
\label{ub:diffuse}
\limsup_{n\to\infty} \Eeen(v_n) \leq \Fzb(v_0).
\end{equation}
If in addition $v_0$ minimizes $\Fzb$ under constrained mass and $\e_n \leq C_2(\eta_n)$, then this sequence also satisfies
\begin{equation}
\label{ub:diffuse-nextlevel}
\limsup_{n\to\infty} \Jeen(v_n) \leq \Hzb(v_0).
\end{equation}
\end{itemize} 
\end{theorem}

\bigskip





\bigskip

\begin{remark}\label{slave}
{\bf Choice of the slaving of $\e$ to $\eta$.} 
There are two separate arguments connecting the two parameters:
\begin{itemize}
\item If the sharp-interface approximation is to be reasonable, then the scaling should be such that the interfacial width is small with respect to the size of the particles. Since a particle has diameter $O(\eta)$, this translates into the condition $\e\ll \eta$. 
\item $\Fzb$ is infinite on structures that are \emph{not} collections of point masses. If $\Fzb$ is to be the limit functional of $\Eee$, then along any sequence that {does not} converge to such point-mass structures $\Eee$ should diverge. It turns out that this provides a stronger condition, as we now show. 

For $\Eee$, every function $v$ is admissible. Under constrained mass $M$, an obvious candidate for the limit behavior is the function $v\equiv M$, with  energy scaling $\Eee(1)\sim \eta^{4}/\e$.
On the other hand, if the functional $\Eee$ is close to $\Fzb$, then we will have $\Eee\approx \Fzb = O(1)$.
Therefore the ratio $\eta^{4}/\e$ is critical. If this ratio is small, then the constant state has lower energy than localized states, and we do not expect the functional $\Feb$ to be a good approximation of $\Eee$. On the other hand, if the ratio $\eta^{4}/\e$ is large, then localized states have lower energy than constant states. 
\end{itemize}

In Theorem~\ref{T1} above,  
%
the lower bound is responsible for forcing divergence of the energy along sequences which do not converge to point masses; the lower bound therefore requires $\e\ll \eta^{4}$. The extra factor  $\eta_n^\zeta$ is  used  in the truncation part of proof: in relating a diffuse-interface sequence to a  sharp-interface sequence, we truncate at a suitable level set of the interface,  and the small factor $\eta_n^\zeta$  is used to quantify the closeness in interfacial energies with respect to the surface tension $\sigma$.

For the upper bound, we would  ideally  require $\e_n = o(\eta_n)$. What we assume, $\e_n\leq C_1(\eta_n)$ and $\e_n\leq C_2(\eta_n)$, are stronger requirements. However, at this stage we do not know the exact local behavior for minimizers of  $\fzb$. In two dimensions we can fully characterize this local behavior, and as we shall see in Section \ref{2D}, this allows us to require only the weaker condition  $\e_n = o(\eta_n)$ (up to a logarithmic correction). 
In three dimensions we use a convenient version of the Modica-Mortola profile construction which does not give an optimal scaling in terms of closeness of energies (cf. Lemma \ref{MM-construction}). 
If one can establish the conjectured behavior for the local problem (see Section \ref{localstructure}), one can then achieve a sharper slaving. 
\end{remark}

\section{Previous results for the sharp interface limit}\label{sharp}

In \cite{ChoksiPeletier09} we dealt with the sharp-interface functionals  that arise from letting $\epsilon$ tend to zero for fixed $\eta$.
For  $\Eee$ and $\Jee$ respectively these limit functionals are 
 \begin{equation}
\label{def:Fes}
\Feb \, : = \, \begin{cases}
  \eta\, \sigma  \int_\Tb |\nabla v|  \,\,\, 
+\,\,\, \eta \,  \Hmospace v\Tb ^2 
  &\text{if } v\in BV(\Tb; \{0,1/\eta^3\})\\
  \infty &\text{otherwise}, 
\end{cases}
\end{equation}
and 
\[ 
\Heb(v) := \eta^{-1}\left[ \Feb(v) - \fzb\left(\int_\Tb v\right)\right] .
\]
We  proved  that 
\[
\Feb\, \stackrel{\Gamma}\longrightarrow \, \Fzb  \qquad 
{\rm and} \qquad 
\Heb \, \stackrel{\Gamma}\longrightarrow \, \Hzb, \qquad {\rm as} \,\,\, \eta \rightarrow 0.  
\]
Precisely, 
\begin{theorem}
\label{T2} Let $\eta_n$ be a sequence tending to $0$. 
 \begin{itemize} 
\item (Condition 1 -- the lower bound and compactness) 
Let  $v_n$ be a sequence such that the sequence of energies
$\Fen(v_n)$ is bounded. Then
(up to a subsequence) $v_n\weak v_0$, $\supp v_0$ is countable, and 
\begin{equation}
\label{lb:sharp}
\liminf_{n \to \infty} \Fen(v_n) \geq \Fzb(v_0).  
\end{equation}
If in addition  $\Hen(v_n)$ is bounded, then the limit $v_0$ is a global minimizer of $\Fzb$ under constrained mass, 
$v_0 =  \sum_i  m^i \delta_{x_i}$ where $\{m^i\} \in {\cal M}$ and  
\begin{equation}
\label{lb:sharp-nextlevel}
 \liminf_{n\to\infty} \Hen(v_n) \geq \Hzb(v_0). 
\end{equation}
\item (Condition 2 -- the upper bound) 
Let $\Fzb(v_0)<\infty$ and  $\Hzb(v_0)<\infty$ respectively. 
Then there exists a sequence $v_n \weak v_0$ such that 
\begin{equation}
\label{ub:sharp}
\limsup_{n \to \infty} \Fen(v_n) \leq \Fzb(v_0). 
\end{equation}
If  $\Hzb(v_0)<\infty$ then 
exists a sequence $v_n \weak v_0$ such that 
\begin{equation}
\label{ub:sharp-nextlevel}
 \limsup_{n\to \infty} \Hen(v_n) \leq \Hzb(v_0). 
\end{equation}
\end{itemize} 
\end{theorem}

\bigskip

We  recall from \cite{ChoksiPeletier09} some properties of $\fzb$: 
\begin{enumerate}
\item For every $a>0$, $\fzb'$ is non-negative and bounded from above on $[a,\infty)$. 
\item 
If $\{m^i\}_{i\in\N}$ with $\sum_i m^i<\infty$ satisfies
\begin{equation}
\label{cond:mass-minimizer}
\sum_{i=1}^\infty \fzb(m^i) = \fzb\Bigl(\sum_{i=1}^\infty m^i\Bigr),
\end{equation}
then
 only a finite number of $m^i$ are non-zero.
\end{enumerate}

\begin{remark}\label{remark-limit}
In proving Theorem \ref{T2}, the bulk of the work was confined to the lower-bound inequalities wherein, after establishing compactness, 
one needed a characterization of sequences with bounded energy and mass. The characterization implied that such a sequence eventually consists of  collection of non-overlapping,   well-separated connected components (see~\cite[Lemma 5.2]{ChoksiPeletier09}). 

We  note that in proving the second-order $\Gamma$ convergence, we saw  that   for an admissible sequence  $v_n$,  
the boundedness of  $\Hen(v_n)$ implied both a minimality condition and compactness:
\begin{itemize}
\item The minimality condition arose from the fact that  $\Eeen(v_n)$ must converge to its minimal value and implied that the $\{m^i\}$ must satisfy \pref{cond:mass-minimizer}. Hence by property 2 above,  the number of limiting particles must be finite.  
\item The compactness condition implied that for each $m^i$, the minimization problem defining 
 $\fzb(m^i)$ (namely \pref{def:fe3})  had a solution. 
\end{itemize}
These condition are responsible for the additional properties of the weights $m^i$ (c.f. ${\cal M}$) in the definition of $\Hzb$. 
\end{remark}

\section{Proof of Theorem \ref{T1}}\label{prooftheorem}

The proof of Theorem \ref{T1} relies on  Theorem \ref{T2}. For the lower bound we use a suitable truncation to relate the approximating diffuse-interface sequence to a sharp-interface sequence  with the same limit and whose difference in energy is small. 
For the upper bound, we modify, in a neighborhood of the boundary,   the sharp-interface recovery sequence given by Theorem \ref{T2} via a {\it quantification  of the Modica-Mortola optimal-profile construction} (\cite{ModicaMortola77}). Such a result is provided  by a lemma of Otto and Viehmann   \cite{OttoViehmann09}. 

\begin{lemma} \label{MM-construction}
Let $\alpha >0$. There exists a constant $C_0(\alpha)$ such that for any  characteristic function $\chi$  of a subset of\/ $\T^3$ and   $\delta>0$, there exists 
 an approximation $u \in H^1(\T^n, [0,1])$ with 
 \[ \int_{\T^3} \delta \, |\nabla u|^2 \,\,+ \,\, \frac{1}{\delta} \, u^2\, (1 - u^2) \, dx \,\,\le \,\, (\sigma + \alpha) \int_{\T^3} |\nabla \chi|, \]
 and 
 \[   \int_{\T^3} \left\vert \chi - u \right\vert \, dx \,\, \le \,\, C_0(\alpha) \, \delta \,  \int_{\T^3} |\nabla \chi|. \] 
\end{lemma}
\noindent 

The proof of Lemma \ref{MM-construction}  follows from the proof of 
 Proposition 1,  Section 7 in \cite{OttoViehmann09}. 
Note that  in \cite{OttoViehmann09},   the authors deal with the   functional   
\[
\int_{\Omega} \frac{\delta}{2(1 - u^2)} \, |\nabla u|^2 \,\,+ \,\, \frac{1}{2\delta} (1 - u^2) \, dx,   
\]
 defined    on cubes of arbitrary size $\Omega$.
Here the wells are at $\pm 1$ and more importantly, this scaling produces unity as the limiting surface tension $\sigma$. However the structure of their proof only uses the fact that this functional $\Gamma$-converges to 
\[ \int_{\Omega} |\nabla u|. \] 
Hence our Lemma  \ref{MM-construction} follows directly not from the statement of their Proposition 1 but from its proof.

\bigskip

\begin{proof}[Proof of Theorem \ref{T1}]
We first prove Condition 1 (the compactness and lower bounds). Let $\e_n$, $\eta_n$, and $v_n$ be sequences as in the theorem, such that $\Eeen(v_n)$ is bounded (but not necessarily $\Jeen(v_n)$, yet). For part of the proof we will work with the sequence and the energy in the original scaling  $u_n$, given by $u_n = \eta_n^3 v_n$. In terms of $u_n$, we find 
\[
\Eeen(v_n) = \frac{\e_n}{\eta_n^{2}}\int_{\T^3} |\nabla u_n|^2 
  + \frac1{\eta_n^{2}\e_n} \int_{\T^3} W(u_n) 
  + \frac{1}{\eta_n^5} \Hmo {u_n}^2,
\]

Following~\cite{ModicaMortola77} we define the continuous and strictly increasing function 
\[
\phi(s) := 2\int_0^s \sqrt{W(t)}\, dt,
\]
and note that as a consequence of the inequality  $a^2 + b^2 \ge 2ab$, we have 
\begin{equation}
\label{lowerbound:diffuse1}
\Eeen(v_n) \geq \frac1{\eta_n^{2}} \int_{\T^3} |\nabla \phi(u_n)|
  + \frac{1}{\eta_n^5}\Hmo {u_n}^2.
\end{equation}

Now set $\alpha_n = 1/(\sigma-\eta_n^\zeta)$, where as before $ \sigma  = \ 2 \int_0^1 \sqrt{W(t)} \, dt = \phi(1)-\phi(0)$. Fix  $\delta_n>0$ by the condition 
\[
\phi(1-2\delta_n)-\phi(2\delta_n) = \phi(1)-\phi(0)-\eta_n^\zeta = \frac1{\alpha_n},
\]
and note that the quadratic behavior of $W$ at $0$ and $1$ implies that $\delta_n = O(\eta_n^{\zeta/2})$. We also introduce the notation $[u]$ for the clipping to the interval $[0,1]$:
\[
[u] := \min\{1,\max\{0,u\}\}.
\]

Using the characterization of perimeter (cf.~\cite{FlemingRishel60} or~\cite[Th.~2.1]{Ambrosio89})  as
\[
\int_{\T^3} |\nabla \phi([u_n])| = \int_{\phi(0)}^{\phi(1)} \H^1(\partial^*\{\phi([u_n])>t\})\, dt,
\]
we estimate the size of the set 
\[
A_n := \left\{t\in[\phi(0),\phi(1)]: 
  \H^1(\partial^*\{\phi([u_n])>t\}) \geq \alpha_n \int_{\T^3} |\nabla \phi([u_n])|\right\}
\]
by
\[
|A_n| = \int_{A_n} 1\, dt 
  \leq \frac1{\alpha_n\int_{\T^3} |\nabla \phi([u_n])|} 
    \int_{\phi(0)}^{\phi(1)} \H^1(\partial^*\{\phi([u_n])>t\}) \, dt  = \frac1{\alpha_n}.
\]
By the definition of $\alpha_n$ and $\delta_n$ it follows that there exists a $t_n\in [\phi(\delta_n), \phi(1-\delta_n)]\setminus A_n$, for which therefore
\begin{equation}
\label{est:levelcurves}
\H^1(\partial^*\{\phi([u_n])>t_n\}) < \alpha_n \int_{\T^3} |\nabla \phi([u_n])|
\leq \alpha_n \int_{\T^3} |\nabla \phi(u_n)|.
\end{equation}

We now construct an auxiliary sequence $\overline u_n$ such that the corresponding $\overline v_n = \overline u_n/\eta_n$ will be admissible for the sharp-interface functional $\Feb$. We map the values of $u_n$ to $\{0,1\}$ with cutoff $\phi^{-1}(t_n)$:
\[
\overline u_n(x) := \begin{cases}
0 & \text{if }\phi(u_n(x))<t_n\\
1 & \text{if }\phi(u_n(x))\geq t_n,
\end{cases}
\]
so that
\begin{equation}
\label{eq:levelcurves2}
\int |\nabla \overline u_n| = \H^1(\partial^*\{\phi([u_n])>t_n\}).
\end{equation}
We estimate the difference in $L^2$ and $H^{-1}$ of $u_n$ and $\overline u_n$.  
Since $\phi^{-1}(t_n)\in[\delta_n,1-\delta_n]$, the function
\[
\psi_n(u) := \begin{cases}
\;u^2 & \text{if } \phi(u) < t_n\\
(1-u)^2 & \text{if }\phi(u)\geq t_n
\end{cases}
\]
is bounded from above by an increasing factor times $W$,
\[
\psi_n(u) \leq C\delta_n^{-2} W(u) \leq C'\eta_n^{-\zeta}W(u)\qquad
\text{for some $C$, $C'$ independent of $n$}.
\]
Therefore the sequences $u_n$ and $\overline u_n$ are close in $L^2$:
\[
\|u_n-\overline u_n\|_{L^2}^2 = \int_{\T^3} \psi_n(u_n)
\leq C'\eta_n^{-\zeta} \int_{\T^3} W(u_n) = O(\e_n\eta_n^{2-\zeta})\to0,
\]
where the final estimate results from the boundedness of $\Eeen(v_n)$. Consequently they are also close in $H^{-1}$, 
\def\td{\textstyle\dashint}
\begin{eqnarray}
\label{est:Hmo-L2}
\|u_n-\overline u_n - \td(u_n-\overline u_n)\|_{H^{-1}}
& \leq & 
C\|u_n-\overline u_n - \td(u_n-\overline u_n)\|_{L^2}\nonumber\\
& \leq&  C\|u_n-\overline u_n \|_{L^2} \nonumber\\
& = & O(\e_n^{1/2}\eta_n^{1-\zeta/2})\to 0,  
\end{eqnarray}
and the same holds for the squared norms:
\begin{align}\label{est:Hmo-L-3}
&\left|\Hmo{u_n}^2 - \Hmo{\overline u_n}^2\right| \nonumber \\
&\qquad\leq 
  \left(\Hmo{u_n} + \Hmo{\overline u_n}\right)\|u_n-\overline u_n - \td(u_n-\overline u_n)\|_{H^{-1}} \nonumber \\
&\qquad\leq 
\left(2\Hmo{u_n} + \|u_n-\overline u_n - \td(u_n-\overline u_n)\|_{H^{-1}}\right)
   O(\e_n^{1/2}\eta_n^{1-\zeta/2}) \nonumber \\
&\qquad = \left(\eta_n^5 \, \Eeen(v_n)\right)^{1/2}\, O(\e_n^{1/2}\eta_n^{1-\zeta/2})
  + O(\e_n\eta_n^{2 - \zeta}) \nonumber  \\
&\qquad = O(\e_n^{1/2}\eta_n^{7/2-\zeta/2})
  + O(\e_n\eta_n^{2-\zeta})\nonumber  \\
&\qquad = o(\eta_n^{6}),
\end{align}
Note that in the last  lines of \pref{est:Hmo-L2} and \pref{est:Hmo-L-3}, we have used  the hypothesis $\e_n = o(\eta_n^{4 + \zeta})$.

Using~\pref{est:levelcurves} and~\pref{eq:levelcurves2} we transfer the lower bound~\pref{lowerbound:diffuse1} to the sequence $\overline u_n$:
\begin{eqnarray}\label{lbcomparison}
\Eeen(v_n) 
&\stackrel{\eqref{lowerbound:diffuse1},\pref{est:levelcurves}}{\geq}& \frac1{\alpha_n\eta_n^{2}} 
\H^1(\partial^*\{\phi([u_n])>t_n\})
 + \frac{1}{\eta_n^5} \Hmo{u_n}^2\nonumber\\
  & \stackrel{\pref{eq:levelcurves2}, \pref{est:Hmo-L-3}}{=}& \frac1{\alpha_n\eta_n^{2}} \int_{\T^3} |\nabla \overline u_n|
  + \frac{1}{\eta_n^5} \Hmo{\overline u_n}^2
  + o(\eta_n)\nonumber\\
&=& \frac{\eta_n}{\alpha_n} \int_{\T^3} |\nabla \overline v_n|
  + {\eta_n}\Hmo{\overline v_n}^2 
  + o(\eta_n)\nonumber\\
&\geq& \frac1{\sigma \alpha_n} \Fen(\overline v_n) + o(\eta_n), 
\end{eqnarray}
where in the last line we used the fact that $\sigma \alpha_n >1$ (note that  $\sigma \alpha_n \rightarrow 1$ as $n \rightarrow \infty$). 

From~\eqref{lbcomparison} it follows that the sequence $\overline v_n$ satisfies the conditions of Theorem~\ref{T2}. Therefore there exists a subsequence $\overline v_{n_k}$ converging to a limit $v_0$, with countable support, such that 
\begin{equation}
\label{lb:transition}
\liminf_{k\to\infty} \sfE_{\eta_{n_k}}(\overline v_{n_k}) \geq \Fzb(v_0).
\end{equation}
The corresponding subsequence $v_{n_k}$ of the sequence $v_n$ also converges weakly to the same limit, since for $\varphi\in C({\T^3})$, 
\[
\left|\int_{\T^3} (v_{n_k}-\overline v_{n_k})\phi\right|
 \leq \frac1{\eta_{n_k}^3}\|u_{n_k} - \overline u_{n_k}\|_{L^2} \|\varphi\|_{L^2}
 = O(\e_{n_k}^{1/2}\eta_{n_k}^{-2 - \zeta/2}) \to 0.
\]
This proves the compactness of the sequence $v_n$ and the characterization of the support of the limit $v_0$. The lower-bound inequality~\eqref{lb:diffuse} then follows from~\eqref{lbcomparison} and~\eqref{lb:transition}.

\bigskip
We address the lower bound for $\Jee$.  We note that boundedness of $\Jeen(v_n)$ implies boundedness of $\Eeen(v_n)$, so that the characterization of the convergence of the sequence given above applies. In addition, by \pref{lbcomparison}, we have 
\begin{align*}
\Jeen(v_n) &= \frac{1}{\eta_n}\left[\Eeen(v_n) - \lscfz\left(\int_{\T^3} v_n\right)\right] \\
&\geq \frac1{\eta_n}\left[\Fen(\overline v_n) - \lscfz\left(\int_{\T^3} v_n\right)\right]  
 + \frac1{\eta_n}\left(\frac1{\sigma \alpha_n}-1\right) \Fen(\overline v_n) 
 + o(1).
\end{align*}
Since $\sigma \alpha_n = 1+o(\eta_n^\zeta)$, with $\zeta>1$, the lower bound \pref{lb:sharp-nextlevel} for $\Heb$ implies 
\[
\liminf_{n\to\infty} \Jeen(v_n) \geq \Hzb(v_0),
\]
which is~\pref{lb:diffuse-nextlevel}.

\bigskip

\bigskip

We now turn to the upper bound (condition 2),  treating $\Eee$ first. As in the proof of Theorem~\ref{T2},   it is sufficient to prove that for any $v_0$ of the form
\[
v_0 = \sum_{i=1}^N m^i\delta_{x^i}, \qquad \text{with $x^i$ distinct},
\]
there exists a sequence $\overline v_n\weak v_0$ with
\begin{equation}\label{ub-goal} 
\limsup_{n\to\infty} \Eeen(\overline v_n) \,\leq \, \Fzb(v_0). 
\end{equation}
See \cite{ChoksiPeletier09} for an explanation. 
Given such a  $v_0$, Theorem~\ref{T2} (specifically \pref{ub:sharp}) provides an admissible  sequence $v_n \weak v_0$ for 
$\Feb$ with 
\begin{equation}\label{ub-sharp} 
\lim_{n\to\infty} \Fen(v_n) = \Fzb(v_0).
\end{equation}
We write $u_n  := \eta_n^3 v_n$, 
which is  the characteristic function of a subset of $\T^3$ composed of  
$N$ sets whose diameters are decreasing to zero. 
For each $n$, Lemma \ref{MM-construction} with $\alpha = \eta_n$ implies that there exists a $C_0(\eta_n)$ such that 
 for any $\e_n>0$, we have an approximation $\overline u_n \in H^1 (\T^3, [0,1])$ such that 
 \begin{equation}\label{ub-MM} 
 \int_{\T^3} \e_n \, |\nabla \overline u_n|^2 \,\,+ \,\, \frac{1}{\e_n} \, \overline u_n^2 \, (1 - \overline u_n^2) \, dx \,\,\le \,\, (\sigma + \eta_n) \int_{\T^3} |\nabla u_n|, 
 \end{equation}
 and 
 \[   \int_{\T^3} \left\vert \overline u_n - u_n \right\vert \, dx \,\, \le \,\, C_0(\eta_n) \, \e_n \,  \int_{\T^3} |\nabla u_n|. \] 
Now let 
\[ \overline v_n \, = \, \frac{\overline u_n}{\eta_n^3}. \]
We have 
\begin{eqnarray}\label{L^1-close}
\| \overline v_n - v_n \|_{L^1(\T^3)} & = & \frac{1}{\eta_n^3}   \int_{\T^3} \left\vert \overline u_n - u_n \right\vert \, dx  \nonumber \\
& \le & \frac{C_0(\eta_n) \e_n}{\eta_n^3} \,  \int_{\T^3} |\nabla u_n| \nonumber \\
& \le & C\,  \frac{C_0(\eta_n) \e_n}{\eta_n }.
\end{eqnarray}
We  will slave $\e_n$ to $\eta_n$ such that the above tends to zero as $n$ tends to infinity. In particular,  $\overline v_n$ and $v_n$ will have the same limit $v_0$. 
We crudely estimate the $H^{-1}$-norm as follows  
\begin{eqnarray}\label{slave-2}
 \left\| v_n - \overline v_n - \textstyle\dashint (v_n - \overline v_n) \right\|^2_{H^{-1}(\T^3)}
& \le & C \left\| v_n - \overline v_n -\textstyle \dashint (v_n - \overline v_n) \right\|^2_{L^2(\T^3)} \nonumber \\
& \le & C \,  \left\| v_n - \overline v_n  \right\|^2_{L^2(\T^3)} \nonumber \\
& \le & C \,  \| v_n - \overline v_n \|_{L^\infty(\T^3)} \,\, \| v_n - \overline v_n \|_{L^1(\T^3)} 
\nonumber \\
& \stackrel{\pref{L^1-close}}{\le} & C \ \frac{C_0(\eta_n) \e_n}{\eta^4_n }. 
\end{eqnarray}
Next we note that 
\begin{eqnarray}\label{ub-main}
\Eeen(\overline v_n)  & = & \frac{\e_n}{\eta_n^{2}}\int_{\T^3} |\nabla \overline u_n|^2 
  + \frac1{\eta_n^{2}\e_n} \int_{\T^3} W(\overline u_n) 
  + \frac{1}{\eta_n^5} \Hmo {\overline u_n}^2  \nonumber \\
  & = & 
\frac{1}{\eta_n^2} \,  \int_{\T^3} \Bigl(\e_n \, |\nabla \overline u_n|^2 \,\,+ \,\, \frac{1}{\e_n} \overline u^2(1 - \overline u_n^2)\Bigr) \, dx 
 \quad + \quad \eta_n    \Hmo {\overline v_n}^2  \nonumber\\
& \le & 
\frac{1}{\eta_n^2} \, \int_{\T^3}\Bigl( \e_n \, |\nabla \overline u_n|^2 \,\,+ \,\, \frac{1}{\e_n} \overline u^2(1 - \overline u_n^2) \Bigr)\, dx \quad + \quad \eta_n  \Hmo { v_n}^2  \quad  \nonumber \\
& & \qquad \qquad \qquad\qquad \qquad+ \quad 
\eta_n \, \left\| v_n - \overline v_n - \textstyle\dashint (v_n - \overline v_n) \right\|^2_{H^{-1}(\T^3)} \nonumber \\
& \stackrel{\pref{ub-MM}, \pref{slave-2}}{\le} & \eta_n \, (\sigma \, + \, \eta_n) \int_{\T^3} |\nabla v_n| \,\, + \,\,  \eta_n  \Hmo { v_n}^2  \,\, + \,\, 
C \ \frac{C_0(\eta_n) \e_n}{\eta^3_n }  \nonumber\\
& = & \Feb (v_n) \,\, + \,\, \eta_n^2  \int_{\T^3} |\nabla v_n|  \,\,+\,\, C \ \frac{C_0(\eta_n) \e_n}{\eta^3_n }.
\end{eqnarray} 
Thus  we assume 
\begin{equation}\label{slave-3}
\frac{C_0(\eta_n) \e_n}{\eta^3_n } \,\,\rightarrow \,\, 0 \quad \hbox{\rm as } \quad  n \,\rightarrow  \, \infty,
\end{equation}
and we choose a function $C_1$ as in the Theorem such that~\eqref{slave-3} is satisfied whenever $\e_n\leq C_1(\eta_n)$.
 We now take the limsup as $n \rightarrow \infty$ in \pref{ub-main}, and hence \pref{ub-sharp} gives \pref{ub-goal}. 

\bigskip

For the next order, let 
\[ v_0 \, = \, \sum_{i = 1}^{N} m^i \delta_{x^i}, \qquad \{m^i\} \in {\cal M}. \]
Theorem  \ref{T2} (specifically \pref{ub:sharp-nextlevel})  gives a sequence $v_n \weak v_0$ such that 
\begin{equation}\label{ub-nextsharp}
\lim_{n \rightarrow \infty}   {\rm F}_{\eta_n}  (v_n) \, =  \, \Hzb (v_0). 
\end{equation}
We take $\overline v_n$ to be the diffuse-interface approximation used in the previous upper-bound argument but now taking $\alpha$ to be $\eta_n^2$. 
Hence $\overline v_n \weak v_0$ and  following the steps of \pref{ub-main}, we have 
\begin{equation}\label{ub-nextsharp2}
\Eeen(\overline v_n)  \, \le \,  \Feb (v_n) \,\, + \,\, \eta_n^3  \int_{\T^3} |\nabla v_n|  \,\,+\,\, C \ \frac{C_0(\eta_n^2) \e_n}{\eta^3_n },
\end{equation}
and
\begin{eqnarray*}
\Jeen(\overline v_n) & = & \eta_n^{-1}\left[ \Eee(\overline v_n) - \fzb\left(\int_\Tb \overline v_n\right)\right]\\
& \stackrel{\pref{ub-nextsharp2}}{\le} & 
 \eta_n^{-1} \,\, \left[  \Feb (v_n) \,\, + \,\, \eta_n^3  \int_{\T^3} |\nabla v_n|  \,\,+\,\, C \ \frac{C_0(\eta_n^2) \e_n}{\eta^3_n }  \right.\\
 & & \quad \left. \,\, - \,\, 
\fzb\left(\int_\Tb  v_n\right) \,\, + \,\, \left( \fzb\left(\int_\Tb  v_n\right)  \, - \, \fzb\left(\int_\Tb \overline  v_n\right) \right) \right] \\
& \le & 
{\rm F}_{\eta_n} (v_n) \,\, + O(\eta_n) \,\, + \,\, \eta_n^{-1}\left[ L \, \|\overline v_n - v_n\|_{L^1}  \, + \, C \ \frac{C_0(\eta_n^2) \e_n}{\eta^3_n } \right],
\end{eqnarray*}
where $L$ is the local Lipschitz constant of  $\fzb$ (cf. Remark \ref{remark-limit}). 
Thus  choosing  $\e_n$ such that 
\begin{equation}\label{ub-slaving2}
\frac{C_0(\eta_n^2) \e_n}{\eta^3_n }  \,  \rightarrow \, 0 \quad {\rm as} \quad n \rightarrow \infty, 
\end{equation} 
 \pref{ub-nextsharp} implies 
\[ \limsup_{n\rightarrow \infty} \Jeen (\overline v_n) \,\, \le \,  \, \Hzb (v_0). \]
We choose a function $C_2$ as in the Theorem such that $\e_n\leq C_2(\eta_n)$ implies~\eqref{ub-slaving2}.
\end{proof}

\section{The local structure of minimizers and the variational problem that defines $\fzb$} 
\label{localstructure}

Simulations of minimizers of the diblock copolymer problem show phase boundaries which resemble 
constant mean curvature surfaces (see for example \cite{ChoksiPeletierWilliams09} and the references therein): 
In the regime of this article,  we observe  spherical boundaries. 
Experimental observations in diblock copolymer melts also support this~\cite{Thomasetal88}.  
 On the other hand one can see, for example via vanishing first variation, that on a finite domain the nonlocal term  
will have an effect on the structure of the phase boundary~\cite{ChoksiSternberg06}. While rigorous results on this effect remain open, one would expect  that 
exploiting a small parameter might prove useful and,  
indeed, this is exactly what our first order asymptotics have done: in proving the first order lower bound, we have 
reduced the local optimal shape of the particles to solutions of the variational problem~\eqref{def:fe3} that defines~$\fzb$. 
The details of this calculations can be found in \cite{ChoksiPeletier09}. 
Let us now comment of this problem and present some conjectures. 

We briefly recall the problem defining  $\fzb$.  For $m>0$, minimize 
\[  \int_{{\Bbb R}^3} |\nabla u| \,\, + \,\, \int_{{\Bbb R}^3}  \int_{{\Bbb R}^3}  \frac{u(x) \, u(y)}{4 \pi |x - y|} \, dx \, dy 
\qquad \hbox{\rm over all } \, u \in BV({\Bbb R}^3, \{0,1\}) \,\,\,  {\rm with } \, \,\, \int_{{\Bbb R}^3} u \, dx \, = \, m. \]
Note that the two terms are in direct competition: balls are {\it best} for the first term and {\it worst} for the second\footnote{  
The latter point has an interesting history.~{ Poincar\'{e}} \cite{Poincare1887, Poincare1902}  considered the problem of determining possible shapes of a fluid body of mass $m$ in equilibrium. Assuming vanishing total angular momentum, the total potential energy in terms of $u$, the characteristic function of the body,  is given by 
\[ ({\cal P})\qquad   \int_{{\Bbb R}^3}  \int_{{\Bbb R}^3}  -  \frac{u(x) \, u(y)}{ C \, |x - y|} \, dx \, dy,   \]
where $-(C |x - y|)^{-1}$, $C>0$ is the potential resulting from the gravitational attraction between two points $x$ and $y$ in the fluid. 
Poincar\'{e} showed under some smoothness assumptions that a body has the lowest energy if and only if it is a ball. 
He referred to  some previous work of { Liapunoff} but was critical of its incompleteness. It was not until almost a century later that the essential details were sorted out wherein the 
 heart of  proving the statement  lies in rearrangement ideas of {Steiner} for the isoperimetric inequality. These ideas are 
captured in  the {\it Riesz Rearrangement Inequality} and its development (c.f. \cite{LiebLoss01}): for functions $f,g$ and $h$ defined on ${\Bbb R}^d$, 
\[ \int_{{\Bbb R}^d}  \int_{{\Bbb R}^d}  f(y)\,g(x - y)\,h(x) \, dy \, dx \,\, \le \,\,  \int_{{\Bbb R}^d}  
\int_{{\Bbb R}^d}  f^\ast(y)\,g^\ast(x - y)\,h^\ast(x) \, dy \, dx, \]
where $f^\ast, g^\ast, h^\ast$ denote the spherically decreasing rearrangements of $f,g,h$. While the general case of equality was treated by Burchard in \cite{Burchard96}, for the problem at hand where the function $g \sim |\cdot|^{-1}$ is fixed and symmetrically decreasing, the inequality with the specific case of equality was treated by Lieb in \cite{Lieb77},  thus proving that  
balls are the unique minimizers for the potential problem (${\cal P}$). }.
The function $\fzb(m)$ denotes this minimal value, i.e. 
\[ 
\fzb (m) \, : = \, \inf  \left\{ \int_{{\Bbb R}^3} |\nabla u| \, + \, \int_{{\Bbb R}^3}  \int_{{\Bbb R}^3}  \frac{u(x) \, u(y)}{4 \pi |x - y|} \, dx \, dy
 \,\, \bigg\vert \, \, 
u \in BV({\Bbb R}^3, \{0,1\}), \int_{{\Bbb R}^3} u \, dx \, = \, m \right\}. 
\]
We also define the energy of one ball of volume $m$: 
\[ f(m) : =  (36\pi)^{1/3}m^{2/3} \, + \, \frac{2}{5}\left(\frac{3}{4\pi}\right)^{2/3}  m^{5/3}. \]
Clearly, we  have 
$\fzb (m) \, \le \, f(m).$
We conjecture  the following  scenario. There exists $m^\ast > 0$, such that for all $m \le m^\ast$, there exists a global minimizer associated with $\fzb (m)$, and it is a single ball of mass $m$. For $m > m^\ast$, a minimizer fails to exist. In fact, as $m$ increases past  $m^\ast$, the ball remains a local minimizer, but a minimizing sequence consisting of two balls of equal size that move away from each other has lower limiting energy.  This separation is driven by the $H^{-1}$ interaction energy, which attaches a positive penalty to any two objects at finite distance from each other.  The limiting energy of such a sequence is simply 
the sum of the energies of two non-interacting balls, i.e. $2 f(m / 2)$. The critical $m^\ast$ then is the only positive zero of $f(m) - 2f(m/2)$, $m^\ast \approx 22.066$.

As $ m$ further increases above a certain $m^{\ast\ast}>m^\ast$, a sequence consisting of three balls of equal size is a minimizing sequence for $\fzb(m)$, with limiting value $3f(m/3)$; and so on for higher values of $n$. Specifically,  we conjecture the following: 

\begin{conjecture}\label{conjecture1}
The minimizer associated with  $ \fzb (m)$ exists  iff $m \le m^\ast$, and it is a ball of mass $m$. 
Moreover, for all $m >0$, we have 
\[  
\fzb (m) \, = \inf_{n\in \N} nf(m/n).
\]
The infimum is achieved iff $m\leq m^\ast$. 
\end{conjecture}

One might ask as to what {\it is} known  about global minimizers. In our  previous article  \cite{ChoksiPeletier09} on the sharp-interface functionals,  we prove that if a sequence (in $\eta$) has bounded energy  $\Heb$ then it must converge to a weighted sum of delta functions where all the weights $m^i$ must have a corresponding  minimizer of  $\fzb(m^i)$.  One can readily check, via trial functions, that such a sequence exists. 
Thus for certain values of $m$, a minimizer of $\fzb (m)$ does exist. 
 Unfortunately, our lower bound compactness  argument  gives no explicit range for the possible limiting  weights $m^i$.  
One could also consider local minimizers, and in particular one can study the stability of balls. A calculation of the second variation using spherical harmonics  (unpublished) indicates that the ball retains stability up to $m_c \approx 62.83$, 
well past the critical mass $m^\ast$.  

Proving Conjecture \ref{conjecture1} would for the first time provide some rigorous justification for why minimizers of the diblock copolymer problem have phase boundaries which resemble periodic constant mean curvature surfaces, supporting the idea that at small length scales, the perimeter (short-range) effects override  the nonlocal (long-range) effects.

\section{Analogous results in two dimensions}\label{2D}
As in \cite{ChoksiPeletier09}, we summarize the analogous results for $d = 2$. While we do not give all the details, we give the essential features which should enable the reader to complete the proofs.  
The fundamental difference between two and three dimensions is that the $H^{-1}$-norm is {\it critical} in two dimensions. 
As explained in \cite{ChoksiPeletier09}, after rescaling with $v = u / \eta^2$ this involves slaving $\gamma$ to~$\eta$ via 
\[ \gamma = \frac{1}{\logeta \eta^3}, \]  
and the two-dimensional function analogous to $\Eee$ becomes  
\[ \Eeet (v) :=  {\e\eta^{3}}\int |\nabla v|^2
             + \frac{\eta^{3}}\e \int \widetilde W(v) 
             +\invlogeta\Hmo v^2.
\]
Here the rescaled double-well energy is now  
\[
\widetilde W(v) := v^2 (1-\eta^2 v)^2.
\]
The analogous  sharp-interface ($\e \rightarrow 0$) limit is given by 
\[
\Fea (v) \, := \, \begin{cases}
\sigma \,  \eta\int_\T |\nabla v| 
+\invlogeta\Hmospace v\T ^2 
  &\text{if } v\in BV(\T, \{0,1/\eta^2\})\\
  \infty &\text{otherwise},
\end{cases}
\]
where $\sigma$ is again given by \pref{st}.
The  first-order limit is defined by 
\[ 
\Fza(v) := \begin{cases}
\sum_{i\in I} \lscfza(m^i) & \text{if } v = \sum_{i\in I} m^i \delta_{x^i} \text{ with } I \text{ countable, }\{x^i\} \text{ distinct, and }m^i\geq0\\
\infty & \text{otherwise}.
\end{cases}
\]
where the function $\lscfza:[0,\infty)\to[0,\infty)$ is defined as follows. 
 Let 
 \begin{align}
\fza(m) &:= \frac{m^2}{4\pi} + 2\sigma \sqrt{\pi m}\notag\\
&= \frac{m^2}{4\pi} + \inf\left\{\sigma\int_{\R^2} |\nabla z|: 
   z\in BV(\R^2;\{0,1\}), \int_{\R^2} z = m \right\}.
\label{def:fe}
\end{align}
An interesting feature here is the explicit nature of $\fza$ (in contrast to \pref{def:fe3}).  
The first term is the dominant part of the $H^{-1}$ norm in two dimensions, and it arises from the fact that the logarithm is {\it additive} with respect to multiplicative scaling. We introduce the lower-semicontinuous envelope function
\begin{equation}
\label{def:lscfz}
\lscfza(m) := \inf \left\{\sum_{j\in J} \fza(m^j): m^j>0, \sum_{j\in J} m^j = m
\right\}.
\end{equation}
 
For the next order, note that 
 \[
\min \left\{\Fza(v): \int_{\Ta} v = M\right\} = \lscfza(M).
\]
We hence recover the next term in the expansion as the limit of $\Fea - \lscfza$, appropriately rescaled, that is of the functional
\[ \Jeet (v) : =  \left|\log\eta\right| \left[\Eeet(v) - \lscfza\left(\int_\Ta v\right)\right]. \] 
Note that the corresponding sharp interface function is 
\[
\Hea(v) := \left|\log\eta\right| \left[\Fea(v) - \lscfza\left(\int_\Ta v\right)\right].
\]

In order to define the  second-order limit, we require some preliminary definitions. 
We first recall a lemma whose proof was presented in \cite{ChoksiPeletier09}. 
\begin{lemma}
\label{lemma:charmin}
Let $\{m^i\}_{i\in \N}$ be a solution of the minimization problem
\begin{equation}
\label{pb:equalmass}
\min \left\{\sum_{i=1}^\infty \fza(m^i): m^i\geq 0, \ 
   \sum_{i=1}^\infty m^i = M.\right\}.
\end{equation}
Then only a finite number of the terms  $m^i$ are non-zero and all the non-zero terms are equal. In addition, if one $m^i$ is less than $2^{-2/3}\pi$, then it is the only non-zero term.  
\end{lemma}
Let 
\[
 f_0(m) := \frac{m^2}{8\pi} \left(3-2\log\frac m\pi\right).
\]
For $n\in \N$ and $m>0$ the sequence $n\otimes m$ is defined by
\[
(n\otimes m)^i := \begin{cases}
  m & 1\leq i\leq n\\
  0 & n+1 \leq i < \infty.
\end{cases}
\]
Let $\widetilde{\mathcal M}$ be the set of optimal sequences for the problem~\pref{pb:equalmass}:
\[
\widetilde{\mathcal M} := \left\{ n\otimes m : n\otimes m \text{ minimizes }\pref{pb:equalmass} \text{ for }M= nm, \text{ and }\lscfza(m)=\fza(m)\right\}.
\]
Then define
\begin{equation}
\label{def:Hz}
\Hza(v): = \begin{cases}
\displaystyle n \left\{ f_0(m) \, +\,  m^2\, g^{(2)} (0) \right\}\,  +& \\
  \qquad {} \,\,  
  \displaystyle\frac {m^2}2  \sum_{\substack{i,j\geq 1\\ i\not=j}} G_\Ta(x^i-x^j)
     &\displaystyle\text{if } v = m\sum_{i=1}^n \delta_{x^i}, \,\, \{x^i\} \text{ distinct}, \, \, n\otimes m \in \widetilde{\mathcal M},\\
\infty & \text{otherwise}, 
\end{cases}
\end{equation}
where the function $g^{(2)}$ was defined in~\pref{eq:G_T-g2}. 
We briefly comment on these functionals and their properties. 
 As in three dimensions, the boundedness of  $\Jeet$ implies that the limiting weights $m^i$ satisfy both a minimality condition and a compactness condition. The minimality condition implies that  $n\otimes m$ minimizes \pref{pb:equalmass}. 
The compactness condition implies that 
\begin{equation}
\label{cond:compactness-2d}
\lscfza(m^i) = \fza(m^i). 
\end{equation}
As we can see from Lemma  \ref{lemma:charmin}, 
 the minimality condition provides a characterization that is stronger than in three dimensions: in particular the masses must be equal. 
Let us also comment on the function $f_0$. 
The minimization problem \pref{def:fe} has only balls (here circular disks) as solutions. Thus in computing the small-$\eta$ asymptotics of 
 $\Jeet$, the $H^{-1}(\R^2)$-norm of a two-dimensional disc of mass $m$ enters. The functional $f_0 (m)$ is exactly this value.

 \begin{theorem}
\label{T3}
\begin{itemize} 
\item (Condition 1 -- the lower bound and compactness) 
Let $\e_n$ and $\eta_n$ be sequences tending to zero such that $\e_n\eta_n^{-3-\zeta}\to0$ for some $\zeta>0$. Let $v_n$ be a sequence such that the energy 
$\Eeent(v_n)$ is bounded. Then
(up to a subsequence) $v_n\weak v_0$, $\supp v_0$ is countable, and 
\begin{equation}
\label{lb:diffuse-2d}
\liminf_{n\to\infty} \Eeent(v_n) \geq \Fza(v_0).
\end{equation}
If in addition $\Jeent(v_n)$ is bounded, then the limit $v_0$ is a global minimizer of $\Fza$ under constrained mass, and 
\begin{equation}
\label{lb:diffuse-nextlevel-2d}
\liminf_{n\to\infty} \Jeent(v_n) \geq \Hza(v_0).
\end{equation}
\item (Condition 2 -- the upper bound)
Let $\e_n$ and $\eta_n$ be sequences tending to zero such that $\e_n\eta_n^{-1}\left|\log\eta_n\right|\to0$. Let $v$ be such that $\Fza(v)<\infty$.
Then there exists a sequence $v_n \weak v$ such that 
\[
\limsup_{n\to\infty} \Eeent(v_n) \leq \Fza(v).
\]
If in addition $v$ minimizes $\Fza$ under constrained mass, and if $\e_n\eta_n^{-1}\left|\log\eta_n\right|^2\to0$, then this sequence also satisfies
\[
\limsup_{n\to\infty} \Jeent(v_n) \leq \Hza(v).
\]
\end{itemize} 
\end{theorem}
 
The proof of Theorem \ref{T3} is very similar to that of Theorem \ref{T1}.  
 Again, we rely heavily on the lower bound estimate and upper bound recovery sequence of the associate sharp interface problems.   
 We summarized those results in \cite{ChoksiPeletier09}. 
 The lower bound inequality follows verbatim the three-dimensional case, the differences in dimension reflected  by the exponent 3 as opposed to 4 in the slaving of  $\e_n$ to $\eta_n$.
 
 The main difference comes in the upper bound and this is reflected in the less restrictive slaving of  $\e_n$ to $\eta_n$. In two dimensions, 
 minimizers associated with the first order limit are necessarily circular droplets. This gives an upper-bound recovery sequence of 
 circular droplets (cf. \pref{def:fe}). To regularize the circular boundaries, one can bypass 
Lemma \ref{MM-construction} and  simply use a one-dimensional optimal profile to approximate a Heaviside function.  The advantage here is the explicit dependence of $C_0$ on $\eta_n$.  
For the analogous step to \pref{slave-2}, 
one can  use the following  interpolation inequality corresponding to the `nearly'-embedding of $L^1$ in $H^{-1}$
to relate the $H^{-1}$-norm to the $L^1$-norm: 
\begin{lemma}
\label{th:LlogL-estimate}
Let $f\in L^\infty(\T^2)$ with $\int_{\T^2} f = 0$. Then there exists a constant $C>0$ such that
\[
\|f\|_{H^{-1}(\T^2)}^2 \leq C\|f\|_{L^1(\T^2)}^2 
  \left(1+\log \frac{\|f\|_{L^\infty(\T^2)}}{\|f\|_{L^1(\T^2)}}\right)
\]
\end{lemma}    

Since the proof this inequality is short and to our knowledge, absent from the literature, we end this section with its proof. To this end, we first derive an inequality proved by Brezis and Merle~\cite{BrezisMerle91} in a slightly different form.
\begin{lemma}
\label{lemma:BrezisMerle}
There exists a constant $C_0\geq 1$ such that 
\[
\int_\Ta e^{|\phi|} \leq C_0,
\]
for all $\phi\in W^{2,1}(\Ta)$ satisfying 
\[
\int_\Ta |\Delta \phi| = 1.
\]
\end{lemma}

\begin{remark}
As the proof below shows, the result holds true for any $\phi$ such that $\int_\Ta |\Delta\phi| < 4\pi $; the constant $C_0$ diverges as the critical value of $4\pi$ is approached.
\end{remark}

\begin{proof}[Proof of Lemma~\ref{lemma:BrezisMerle}]
Setting $f(x) := -\Delta \phi$, so that $\int |f| = 1$,  we have
\[
\phi(x) = \int_\Ta G_\Ta (y) f(x-y)\, dy,
\]
and note that by~\pref{eq:G_T-g2}
\[
|G_\Ta (y) | \leq C - \frac1{2\pi} \log |y|
\]
for some $C>0$, and for all $y\in(-1/2,1/2)^2$. Therefore, using Jensen's inequality, 
\begin{align*}
\int_\Ta e^{|\phi(x)|}\, dx 
  &\leq \int_{(-1/2,1/2)^2} \exp\left(\int_ {(-1/2,1/2)^2} |G_\Ta (y)||f(x-y)|\, dy\right)\, dx \\
  &\leq e^{C} \int_{(-1/2,1/2)^2} \exp\left(\int_ {(-1/2,1/2)^2} \log\left(|y|^{-1/2\pi}\right)|f(x-y)|\, dy\right)\, dx \\
  &\leq e^{C} \int_{(-1/2,1/2)^2} \int_ {(-1/2,1/2)^2} |y|^{-1/2\pi}|f(x-y)|\, dy dx \\
  &= e^{C}  \int_{{(-1/2,1/2)^2}}|y|^{-1/2\pi}\, dx\\
  &=: C_0.
\end{align*}

%
\end{proof}

\begin{proof}[Proof of Lemma~\ref{th:LlogL-estimate}]
Set $\Phi(s) := |s|\log(1+C_0|s|)$ and let $\Phi^*$ be the convex conjugate
$\Phi^*(t) := \sup_{s\in\R} (ts-\Phi(s))$. From the lower bound $\Phi(s)\geq |s|\log(C_0|s|)$ we derive the upper bound 
\[
\Phi^*(t)\leq C_0^{-1}e^{|t|}.
\]
Define the Orlicz norm
\[
\|f\|_\Phi := \inf \left\{\lambda >0: \int_\Ta \Phi\Bigl(\frac f\lambda\Bigr) \leq 1\right\}.
\]
Then we have the H\"older inequality (see, for example, Section 3.3~of \cite{RaoRen73}) 
\[
\int_\Ta fg \leq 2 \|f\|_\Phi \|g\|_{\Phi^*}.
\]

To prove Lemma~\ref{th:LlogL-estimate}
 we take $f\in L^\infty$, $f\not= 0$,  with $\int f =0$, and by multiplying $f$ with a constant we can assume that $\int |f| = 1$. Setting $-\Delta \phi =f$, we have
\[
\|f\|_{H^{-1}(\Ta)} ^2
= \int_\Ta f\phi 
\leq 2 \|f\|_\Phi \|\phi\|_{\Phi^*}
\leq 2 \|f\|_\Phi.\]
The second inequality above follows from remarking that
\begin{align*}
\|\phi\|_{\Phi^*} &= 
\inf \left\{\lambda >0: \int_\Ta \Phi^*\Bigl(\frac \phi\lambda\Bigr) \leq 1\right\}\\
&\leq \inf \left\{\lambda >0: C_0^{-1}\int_\Ta e^{|\phi|/\lambda} \leq 1\right\}\\
&\hspace{-0.5cm}\stackrel{{\rm Lemma} \, \ref{lemma:BrezisMerle}}\leq 1.
\end{align*}
Now let  $\lambda_\ast := \|f\|_\Phi$. Since the map $\lambda \rightarrow \int_\Ta \Phi\Bigl(\frac f\lambda\Bigr)$ is continuous at $\lambda_\ast$, we must have 
\[ \int_\Ta \Phi\Bigl(\frac {f}{\lambda_\ast}\Bigr) \, = \, 1.\]
Thus 
\[
\lambda_\ast = \int_\Ta |f(x)| \log \left(1+C_0\frac{f(x)}{\lambda_\ast}\right)\, dx
\leq \log \left(1+C_0\frac{\|f\|_\infty}{\lambda_\ast}\right)
\]
or
\[
\lambda_\ast (e^{\lambda_\ast} -1 ) \leq C_0\|f\|_\infty.
\]
We note that 
\[
\frac{\log 2}2 \, e^\lambda \leq \lambda(e^\lambda-1)
\qquad\text{for all $\lambda>0$ with }\lambda(e^\lambda-1)\geq 1. 
\]
Hence if $\lambda_\ast (e^{\lambda_\ast} - 1) \geq 1$, then 
\[
\|f\|_{H^{-1}(\Ta)} ^2\leq 2\|f\|_\Phi = 2 \lambda_\ast \leq 2 \log \frac {2C_0\|f\|_\infty}{\log 2},
\]
On the other hand, if $\lambda_\ast (e^{\lambda_\ast} - 1) < 1$, then since 
$\lambda \mapsto \lambda$ $ (e^\lambda -1)$ is increasing,  we have $\lambda_\ast \leq  \bar \lambda$, where $\bar \lambda (e^{\bar\lambda}-1) = 1$.
Since 
\[
\frac{\log 2}2 \, e^{\bar\lambda} \leq {\bar\lambda}(e^{\bar\lambda}-1),   
\]
 we have 
\[
\|f\|_{H^{-1}(\Ta)} ^2\leq 2\|f\|_\Phi = 2\lambda_\ast \leq 2 \bar\lambda \leq 2 \log \frac {2C_0\|f\|_\infty}{\log 2}.
\]
Replacing  $f$ with  $ f /\|f\|_{L^1}$ gives the desired inequality. 

\end{proof}

\section{Discussion, dynamics, and related work}

Together with \cite{ChoksiPeletier09}, we have presented an analysis of the small-volume regime for the {\it diblock copolymer problem}. This has been accomplished by an asymptotic description of the energy functional  in the small volume-fraction regime. 
We refer to the discussion section of \cite{ChoksiPeletier09} for comments on  the role of the mass constraint with respect to the limit functionals and  the fundamental differences between the two- and three-dimensional cases. 
As described above, in three dimensions many open problems remain with respect to the local structure problem and it is here that one should first focus in order to rigorously address the role of the nonlocal term on shape effects. 

This asymptotic study has much in common with the asymptotic analysis of the well-known {\it Ginzburg-Landau functional} for the study of magnetic vortices ({\it cf.} \cite{SandierSerfaty07,  JerrardSoner02, AlbertiBaldoOrlandi05}).  Our problem is much more direct as it pertains to the asymptotics of the support of minimizers. This is in strong contrast to the Ginzburg-Landau functional wherein one is concerned with an intrinsic vorticity quantity which is captured via a certain gauge-invariant Jacobian determinant of the order parameter.   

Our results are consistent with and complementary to two other recent studies in the regime of small volume fraction.  
In \cite{RenWei08} Ren and Wei prove the existence of  sphere-like solutions to the Euler-Lagrange equation of  (\ref{def:Ees}), and further investigate their stability.  They also show that the centers of sphere-like solutions are close to global minimizers of an effective energy defined over delta measures which includes both a local energy defined 
over each point measure, and a Green's function interaction term which sets their location.
While their results are similar in spirit to ours, they are based upon completely different techniques which are local rather than global. 
Very recently. Muratov \cite{Muratov09} proved a strong and rather striking result for the sharp interface problem in two dimensions. 
In an analogous small volume fraction regime, he proves that the global  minimizers  are nearly identical circular droplets of a small size separated by large distances. While this result does not precisely determine the placement of the droplets -- ideally proving periodicity of the ground state, to our knowledge it  presents the first rigorous work characterizing some {\it geometric} properties of the ground state (global minimizer).

We conclude this section on the interesting connection with gradient-flow dynamics. It is convenient to examine either the $H^{-1}$ gradient flow of 
\pref{def:Ees} or the modified Mullins-Sekerka free boundary problem of Nishiura-Ohnishi \cite{NishiuraOhnishi95} which results from taking the gradient flow of the sharp-interface functional. 
In \cite{Helmersetal08, GlasnerChoksi09}  the authors explore the dynamics of small spherical phases (particles).  
By constructing  approximations based 
upon an Ansatz of spherical particles similar to the classical Lifshitz-Slyozov-Wagner theory, one derives a finite dimensional dynamics for particle positions and radii.
Here one finds a separation of  
time scales for the dynamics:  Small particles both exchange material as 
in usual Ostwald ripening, and migrate because of an effectively repulsive nonlocal energetic term.  Coarsening via mass 
diffusion only occurs while particle radii are small, and they 
eventually approach a finite equilibrium size.  Migration, on the other hand, is
responsible for producing self-organized patterns. For large systems, kinetic-type equations which describe the evolution of a probability 
density are constructed.  A separation of time scales between particle growth and 
migration allows for a variational characterization of spatially inhomogeneous quasi-equilibrium states.
Heuristically this matches our findings of (a) a first order energy which is local and 
essentially driven by perimeter reduction, and (b)  a Coulomb-like interaction energy, at the 
next level,  responsible for placement and self organization of the pattern. 
Moreover, in \cite{GlasnerChoksi09}, one finds that both the particle position radii and  centre ODE's have gradient-flow structures  related to energies which can be directly linked  to our first  and second order limit functionals, respectively.  

The natural question is to what extent one can rigorously address the dynamics and the separation of coarsening  and  particle migration effects. Recently, Niethammer and Oshita \cite{NiethammerOshita09} have given a rigorous derivation of the mean-field equations associated with the evolution of radii. Another approach (currently in progress) is via Sandier and Serfaty's  connection between 
$\Gamma$-convergence and an appropriate (weak) convergence of the associated gradient flows 
\cite{SandierSerfaty04, Serfaty09}. Le \cite{ Le09}  has recently used this framework for the $\epsilon \rightarrow 0$ problem,  establishing  convergence of the $H^{-1}$-gradient flow of \pref{def:Ees} to that of the  modified Mullins-Sekerka free boundary problem of Nishiura and Ohnishi \cite{NishiuraOhnishi95}.  While this method gives a rather weak notion of convergence, it allows for much weaker assumptions on the initial data and generic structure of the evolving phases.

\bigskip

{\bf Acknowledgments:} The research of RC  was partially supported by  
an NSERC (Canada) Discovery Grant. The research of MAP was partially supported by 
 NWO project 639.032.306.

\bibliographystyle{acm}
\bibliography{bibliography}

\end{document}